\documentclass[runningheads]{llncs}
\usepackage{amsmath,amsfonts,amssymb}
\usepackage{setspace}
\usepackage{fancyhdr}
\usepackage{lastpage}
\usepackage{extramarks}
\usepackage{chngpage}
\usepackage{soul,color}
\usepackage{graphicx,float,wrapfig}
\usepackage{CJK}
\usepackage{algorithmic}
\usepackage[utf8]{inputenc}
\usepackage{graphicx}
\usepackage{tikz}
\usetikzlibrary{quotes, arrows, graphs}
\begin{document}
\title{A Study on the Splitting Strategy\\ of Server Resources\thanks{Supported by Professor Chenye Wu.}}
%
%
\author{Yiheng Shen}
\authorrunning{Yiheng Shen}
%
\institute{Tsinghua University, IIIS, China \\
\email{shen-yh17@mails.tsinghua.edu.cn}\\
}

\maketitle              
\begin{abstract}
The paper is based on Noar's model of charged queues. We extend this model into multi-server systems with information about length and service rate disclosed for all the customers, and the customers can choose the optimal options. We discuss whether the splitting strategy of the server resource could bring more revenue for the service provider. We prove that any G/D/1 server supplier cannot earn more revenue by splitting his resource under the equal-toll limitations.\\

\end{abstract}
\section{Introduction}
 Suppose that there is a supplier with a certain of computers as server resources. The supplier has two strategies. One is to split its server resources into several parts and make each part works as a server (and finally he owns multiple servers). The other is to construct a single server using all his resources. In real world, we have observed that the first strategy may result in VIP policies. We may consider the company open an extra server for customers with limitations. Since this extra server requires server resource, it is an example for the first strategy. Some other companies simply opens a single server with a very long queue (which is the second strategy). We discuss whether the splitting strategy is beneficial in queueing theory. It is obvious that when we consider the whole throughput, the firstt strategy is better than the second. However, when we consider the revenue and the individual habits of the customers, the case become non-trivial.\\

We first introduce the Naor's model for a single queue. Suppose that there are $n$ classes of customers, each with arrival rate $\lambda_1$, $\lambda_2$,$\cdots$, $\lambda_n$. For each customer in class $i$, let $R_i$ denote the reward he earns after his service is completed and $c_i$ denote the price of his cost per unit of time (second). The size of the jobs for class $i$ is exponentially distributed with rate $1/s_i$ and the arrival for each class is memoriless. The serving rate for the server is $\mu$. Balking is allowed, thus each newly arrived customer chooses to join the queue and after waiting he gets the benefits or to leave the queue without any cost. The toll(admission) price of the server is $\theta$, which means if a customer is entering the queue, he should pay $\theta$ first. The toll is then counted into the total revenue of the server. Each new customer is risk neutral, thus he only take the expected costs into account.

Now we extend this model to a multi-server system, suppose there are $m$ queues in total. Each queue has a server. For a server $j(1\le j\le m)$, we let $\theta_j$ be the toll of this server and $\mu_j$ be the serving rate of this server. Then when a new customer arrives, he may calculate the expected benefits for each queue and he chooses join one of the queues or to leave the system with utility 0. Naturally, each customer is intends to choose the action from which he could get the largest expected utility. Bodas et al.(2011) have also studied the equilibrium in the multi-queue systems with toll price, but in their model, each customer could not see the instantaneous number of customers already in the queue. In our model, each customer could see the length of the queues and then take the individually optimal action.

Since we only care mainly about the revenue received by the supplier of the servers, and the total server resource of the supplier is fixed, we assume that the splitting of the server resource is linear, which means $\sum_{j=1}^m \mu_j= \mu_s$, where $\mu_s$ denotes the total resource of the supplier.

\section{Splitting strategy for G/D/1 queue(s)}
We now prove that we cannot get more benefits from splitting the server resource under some certain limitations.\\

We consider the settings in which the service time for each job is fixed while the arrivals are only independent of the system (not necessarily memoriless), namely \textbf{general fixed-size arrivals}. The toll prices for all the split servers should be identical in splitting strategy. Each server has a toll price and the customers has a waiting cost $c$ and a completion reward $W$. Since the size of each job is fixed, without loss of generality, we set this size to be $1$ and the service rates for the servers should be added up to $\mu$, i.e. the service rate of each job in the OPTn strategy should be $1/\mu$. A customer knows the waiting time in every queue when it arrives at the system, and he will enter the queue with maximum utility or balk with zero utility.

\begin{definition}[equal-toll]
    The property in a splitting strategy that all the split servers have the same toll price.
\end{definition}
\begin{definition}[OPTn strategy]
	The optimal strategy in which we get maximum revenue rate without splitting the server resource.
\end{definition}
\begin{theorem}\label{thm1}  
	No equal-toll splitting strategies outperforms OPTn strategy if the arrivals are general fixed-size arrivals.
\end{theorem}
\begin{lemma}\label{lem1}
	Suppose the system consists of two equal-toll servers $s_1$, $s_2$ and the arrivals to the system are general fixed-size arrivals. Then the throughput of the system cannot exceed a single server whose service rate is $s_1 + s_2$ while toll price is identical to $s_1$ and $s_2$.
\end{lemma}
	We denote the system of the two servers as $S_1$; the system with a single server whose service rate is the sum of the two servers in $S_1$ as $S_2$. We use a list of arrival time $L=[t_1,t_2,\cdots]$ to denote the arrival time of all the customers in the system, namely the arrival schedule list. Suppose the set of the customers which enters the system $S_1$ forms a new schedule list $L'= [t_{i_1}, t_{i_2},\cdots]$. We prove the lemma in two steps:\\

Step 1. If the arrival schedule is exactly $L'$, all the customers in $L'$ could also be accepted by (enter) the system $S_2$. 
\begin{proof}
Suppose service rates of $s_1$ and $s_2$ are $\mu_1$ and $\mu_2$ respectively. The single server at $S_2$ has service rate $\mu_1+\mu_2$. A new customer should wait $f_i(t)$ time if he arrives at $s_i(1\le i\le 2)$ at time point $t$ and should wait $F(t)$ if he arrives at $S_2$. Then we consider the three functions $f_1,\ f_2,\ F$ at the moment right after the time point $t_{i_k}(k\in \mathbb{Z^+})$. It is easy to see that at time point $t=0$, $f_1(0)= \frac{1}{\mu_1}$, $f_2(t)=\frac{1}{\mu_2}$, $F(0)=\frac{1}{\mu_1+\mu_2}$. Let $t_{i_0}^-=t_{i_0}=0$, we can derive the following equations:
\begin{gather*}
\forall k\in \mathbb{Z^+}, \\
f_1(t_{i_k}^-)=\max\{f_1(t_{i_{k-1}})-t_{i_k}+t_{i_{k-1}}, \frac{1}{\mu_1}\}\\
f_1(t_{i_k})=f_1(t_{i_k}^-)+\mathbb{I}[f_1(t_{i_k}^-)\le f_2(t_{i_k}^-)]\cdot \frac{1}{\mu_1},\\
f_2(t_{i_k}^-)=\max\{f_2(t_{i_{k-1}})-t_{i_k}+t_{i_{k-1}}, \frac{1}{\mu_2}\}\\
f_2(t_{i_k})=f_2(t_{i_k}^-)+\mathbb{I}[f_1(t_{i_k}^-)>f_2(t_{i_k}^-)] \cdot \frac{1}{\mu_2}, \\
F(t_{i_k}^-)=\max\{F(t_{i_{k-1}})-t_{i_k}+t_{i_{k-1}}, \frac{1}{\mu_1+\mu_2}\}\\
F(t_{i_k})=F(t_{i_k}^-)+\frac{1}{\mu_1+\mu_2}\cdot \mathbb{I}[\text{The $k^{th}$ customer could enter $S_2$}],
\end{gather*}

where $\mathbb{I}[\text{E}]$ is the indicator function that $\mathbb{I}[\text{E}]=1$ if an event E happens and $\mathbb{I}[\text{E}]=0$ otherwise. $t_{i_k}^-$ denotes the time point just before the time point $t_{i_k}$ (which is the time point that determines which server should the $k^{th}$ customer use). 
\\

We define the following three functions:
\begin{gather*}
    W_1(t) = f_1(t)\cdot \mu_1-1,\\
    W_2(t) = f_2(t)\cdot \mu_2-1,\\
    W(t) = f_2(t)\cdot (\mu_1+\mu_2)-1.
\end{gather*}

These functions represent the sum of the total rest jobs' sizes in $s_1$, $s_2$ and $S_2$, which is also the amount of remaining work in the servers. Furthermore, we get the following relations:
\begin{gather*}
    \forall k\in \mathbb{Z^+},\\
    W_1(t_{i_k}^-)=\max\{W_1(t_{i_{k-1}})-\mu_1\cdot (t_{i_k}-t_{i_{k-1}}),0\},\\
    W_2(t_{i_k}^-)=\max\{W_2(t_{i_{k-1}})-\mu_2\cdot (t_{i_k}-t_{i_{k-1}}),0\},\\
    W(t_{i_k}^-)=\max\{W(t_{i_{k-1}})-(\mu_2+\mu_1)\cdot (t_{i_k}-t_{i_{k-1}}),0\},\\
    W_1(t_{i_k})=W_1(t_{i_k}^-)+\mathbb{I}[f_1(t_{i_k}^-)\le f_2(t_{i_k}^-)],\\
    W_2(t_{i_k})=W_2(t_{i_k}^-)+\mathbb{I}[f_2(t_{i_k}^-)> f_2(t_{i_k}^-)],\\
    W(t_{i_k})=W(t_{i_k}^-)+\mathbb{I}[\text{The $k^{th}$ customer could enter $S_2$}].\\
\end{gather*}

We show that $F(t)\le \min\{f_1(t),f_2(t)\}$ and $W(t)\le W_1(t)+W_2(t)$ at any time point $t= t_{i_k}^-(k\in \mathbb{Z^+})$ by induction.
\begin{itemize}
\item[(i)] When $k=1$, $f_j(t_{i_k}^-)=f_j(0) =\frac{1}{\mu_j}\ (\forall\ 1\le j\le 2)$, $F(t_{i_k})=F(0)= \frac{1}{\mu_1+\mu_2}$. It is obvious that $\frac{1}{\mu_1+\mu_2}<\min\{\frac{1}{\mu_1},\frac{1}{\mu_2}\}$. Since $W(t_{i_k}^-)=W(0)=W_1(0)=W_2(0)=W_1(t_{i_k}^-)=W_1(t_{i_k}^-)=0$, $W(t_{i_k}^-)\le W_1(t_{i_k}^-)+W_2(t_{i_k}^-)$. Therefore the statement holds for $k=1$.
\item[(ii)] Suppose when $k=p\ (p\in \mathbb{Z^+})$, at the time point $t_{i_p}^-$ $(p\in \mathbb{Z^+})$ (just before the $p^{th}$ customer arrives), we have 
\begin{gather*}
F(t_{i_p}^-)\le \min\{f_1(t_{i_p}^-),f_2(t_{i_p}^-)\},\\ W(t_{i_p}^-)\le W_1(t_{i_p}^-)+W_2(t_{i_p}^-)    
\end{gather*}

\newcommand{\tp}{{t_{i_p}^-}}
\newcommand{\tg}{{t_{i_p}}}
as our inductive hypothesis. Then at the time point $t_{i_p}$, since we have the first hypothesis and the waiting cost is identical in all the three servers. By the definition of $L'$, the $p^{th}$ customer could enter the first system, then suppose the reward of this customer is $R_p$ and the waiting cost of this customer is $c_p$ (these two parameters may vary from class to class). $R_p-c_p\cdot\min\{f_1(\tp),f_2(\tp)\}>0$, since $F(\tp)\le \min\{f_1(t_{i_p}^-),f_2(t_{i_p}^-)\}$, we get $R_p-c_p\cdot F(\tp)>0$ thus the $p^{th}$ customer could be accepted by $S_2$ and the indicator function $\mathbb{I}[\text{The $k^{th}$ customer could enter $S_2$}]$ should be 1. Without loss of generality, we suppose that $f_1(\tp)<f_2(\tp)$, then $f_1(\tg) =f_1(\tp)+\frac{1}{\mu_1}$ and $f_2(\tg)=f_2(\tp)$. It is trivial to see that $W(\tg)=W(\tp)+1\le W_1(\tp)+W_2(\tp)+1=W_1(\tg)+W_2(\tg)$. Next we are going to prove $F(\tg)\le \min\{f_1(\tg)+f_2(\tg)\}$. 
\begin{itemize}
    \item [1.] If $f_1(\tg)<f_2(\tg)$, since $F(\tg)=F(\tp)+\frac{1}{\mu_1+\mu_2}<f_1(\tp)+\frac{1}{\mu_1}=f_1(\tg)<f_2(\tg)$, thus $F(\tg)<\min\{f_1(\tg),f_2(\tg)\}$.\\
    \item [2.] If $f_1(\tg)\ge f_2(\tg)$, since $W(\tp)\le W_1(\tp)+W_2(\tp)$, we get 
    \begin{align*}
     F(\tp)\cdot (\mu_1+\mu_2)-1&\le f_1(\tp)\cdot \mu_1 -1+f_2(\tp)\cdot \mu_2-1\\
     F(\tp) + \frac{1}{\mu_1+\mu_2}&\le f_1(\tp)\cdot \frac{\mu_1}{\mu_1+\mu_2}+f_2(\tp)\cdot \frac{\mu_2}{\mu_1+\mu_2}\\
     F(\tp) + \frac{1}{\mu_1+\mu_2}&\le \max\{f_1(\tp),f_2(\tp)\}.
    \end{align*}
    Because $f_1(\tp)<f_2(\tp)$, $F(\tp)+\frac{1}{\mu_1+\mu_2}\le f_2(\tp)$, $F(\tg)=F(\tp)+\frac{1}{\mu_1+\mu_2}\le f_2(\tg)\le f_1(\tg)$ and $F(\tg)\le \min\{f_1(\tg),f_2(\tg)\}$.\\\\
\end{itemize}

Since $\forall a,b\in \mathbb{R}$, $\max\{a+b,0\}\le\max\{a,0\}+\max\{b,0\}$. Let $\Delta t = t_{i_{p+1}}-t_{i_{p}}$, we have
\begin{align*}
 W(t_{i_{p+1}}^-)&=\max\{W(t_{i_{p}})-(\mu_2+\mu_1)\cdot {\Delta t},0\}\\
    &\le \max\{W_1(t_{i_{p}})+W_2(t_{i_{p}})-(\mu_2+\mu_1)\cdot {\Delta t},0\}\\
    &\le \max\{W_1(t_{i_{p}})-\mu_1\cdot {\Delta t},0\}+
    \max\{W_2(t_{i_{p}})-\mu_2\cdot {\Delta t},0\}\\
    &\le W_1(t_{i_{p+1}}^-)+W_2(t_{i_{p+1}}^-).
\end{align*}
Since $F(\tg)-\Delta t \le \min\{f_1(\tg)-\Delta t,f_2(\tg)-\Delta t\}$ and $
  \frac{1}{\mu_1+\mu_2}<\min\{\frac{1}{\mu_1},\frac{1}{\mu_2}\}$, we have 
  \begin{align*}
  F(t_{i_{p+1}}^-)&=\max\{F(t_{i_p})-\Delta t, \frac{1}{\mu_1+\mu_2}\}\\
  &\le \min\{\max\{f_1(\tg)-\Delta t,\frac{1}{\mu_1}\},\ \max\{f_2(\tg)-\Delta t, \frac{1}{\mu_2}\}\}\\
  &= \min\{f_1(t_{i_{p+1}}^-),f_2(t_{i_{p+1}}^-)\}.
  \end{align*}
\end{itemize} 

By induction from (i),(ii), the waiting time at $S_2$ should never surpass the minima of the waiting time in $s_1$ and $s_2$. Thus every job entering the system $S_1$ could enter $S_2$ if the arrival schedule is $L'$.\\
\end{proof}

Step 2. If the arrival schedule is $L$, then the server $S_2$ should have a larger through put compared with situation when the arrival schedule is $L'$. 
\begin{proof}
It can be noticed that the schedule list $L'$ is created after deleting all the customers which balk from $S_1$. Then at any time point in a specific schedule list $L$, the single server $S_2$ could finish more work than the situation when the arrival schedule is $L'$. When the list is $L'$, the throughput of $S_1$ and $S_2$ are the same. However, adding the balked customers into the schedule list $L'$ to change it into $L$ may add throughput to the system $S_2$ while the keep the original throughput of $S_1$. The revenue rate of the supplier is proportional to the throughput of a system when all the servers in the systems are equal-toll, therefore we have proven \textbf{Lemma \ref{lem1}}.
\end{proof}
After proven Lemma 1, we use induction to prove \textbf{Theorem \ref{thm1}}.
\begin{proof}[Theorem \ref{thm1}]
We use induction to prove the theorem. Suppose that in a splitting strategy, the system $S_1$ consists of $m$ servers. 
\begin{itemize}
    \item [(i)] When $m=1$, only one single server has service rate $\mu_s$. Thus the strategy can never do better than OPTn.
    \item [(ii)] Suppose that when $m=p$, splitting the resource into $p$ servers cannot bring more benefit than a single server whose service rate is the sum of the $p$ servers. When there are $m=p+1$ in $S_t$, we now construct a system with two servers $s_1$ and $s_2$, while $s_1$ has a service rate which is the sum of the first $p$ servers' rates and the $s_2$ is same as the $(p+1)^{th}$ server in $S_t$. By induction hypothesis, we get $s_1$ can bring at least the same benefits brought by the first $p$ servers of $S_t$. By \textbf{Lemma \ref{lem1}}, the system with $s_1$ and $s_2$ cannot be better than a single system whose service rate is the sum of $s_1$ and $s_2$. Thus finally we reach the conclusion that $p+1$ servers cannot bring more benefit than a single server with the same total service rate. The statement holds for $m=p+1$.
\end{itemize}

By induction, we reach the conclusion that any equal-toll splitting strategy cannot bring more benefits than a single server under general fixed-size arrivals. Since OPTn is the best strategy for a single server system, we have finished the proof.\\
\end{proof}
Despite of the limitations of \textbf{Theorem \ref{thm1}}, it could be useful in real life when there are very strong limitations on the size of the customers (i.e. in HEYTEA a customer is limited to buy at most two cups of drinks, thus a vast majority of customers would buy two cups of drinks and the service time is nearly fixed to be the time for the workers to make two cups of tea). Since the toll can be set to the average earn per customer, it is also nearly fixed. Therefore, if the shop opens several queues, the toll for each queue is equal and our model would be a good approximation for the system.\\

\section{Conclusion and Future work}
We prove that when the toll for each server is equal and the size of each customer is fixed, the splitting strategy would never bring more benefits, regardless of the class distribution of the customers. 

Since in this paper, we only consider a situation when the size of each job is fixed. It remains open that what is the boundary of gaining more benefit from splitting the server resource. We give a hypothesis that even if the customer(job) sizes are not necessarily fixed, if the toll prices are fixed, the supplier still could not earn revenue from splitting the resources. But in this situation. The step 2 part in the proof of our major lemma would require more complicated work.

How to split the resource to maximize the benefits is also a demanding job. We have only found some examples of gaining more result from some multi-class arrivals, but we have not found out a method to estimate the optimal splitting and to set good toll prices for the split servers. This is mainly because after the splitting, the servers could influence each other and this would make the arrival process at any server hard to describe. 

%

\bibliographystyle{splncs04}
\bibliography{mybibliography}

\end{document}